\input amstex
\documentstyle{amsppt}
\magnification=1200
\hsize=13.8cm
\catcode`\@=11
\def\NoLogo{\let\logo@\empty}
\catcode`\@=\active
\NoLogo

\def\heat{\lf(\frac{\p}{\p t}-\Delta\ri)}
\def \b {\beta}
\def\i{\sqrt{-1}}

\def\lf{\left}
\def\ri{\right}
\def\bbar{\bar \beta}
\def\a{\alpha}
\def\ol{\overline}
\def\g{\gamma}
\def\e{\epsilon}
\def\p{\partial}
\def\delbar{{\bar\delta}}

\def\dbar{\bar\partial}

\def\C{\Bbb C}
\def\R{\Bbb R}

\def\vp{\varphi}

\def\dbar{\bar\partial}
\def\ba{{\bar\alpha}}
\def\bb{{\bar\beta}}
\def\abb{{\alpha\bar\beta}}

\def\i{\sqrt {-1}}

\def \D {\Delta}
\def\aint{\frac{\ \ }{\ \ }{\hskip -0.4cm}\int}
\documentstyle{amsppt}
\vsize=19.0 cm

\leftheadtext{Lei Ni  and Luen-fai Tam}
\rightheadtext{Liouville properties}
\topmatter
\title{Liouville properties of plurisubharmonic functions}\endtitle

\author{Lei Ni\footnotemark and  Luen-Fai Tam\footnotemark}\endauthor
\footnotetext"$^{1}$"{Research partially supported by NSF grant DMS-0203023,
USA.}
\footnotetext"$^2$"{Research partially supported by   Earmarked Grant of
Hong Kong \#CUHK4032/99P.}
\address
Department of Mathematics, University of California, San Diego, La Jolla,
CA 92093
\endaddress
\email{
lni\@math.ucsd.edu}
\endemail

\address
Department of Mathematics, The Chinese University of Hong Kong,
Shatin, Hong Kong, China
\endaddress

\email{lftam\@math.cuhk.edu.hk}
\endemail

\affil
{
University of California, San Diego\\
The Chinese University of Hong Kong
}
\endaffil

\date   October, 2002 (Revised expanded version)
\enddate

\endtopmatter

\document
\subheading{\S0 Introduction}

In this paper we will  prove  a Liouville theorem on
  smooth plurisubharmonic functions on  a complete
noncompact K\"ahler  manifold
 with nonnegative bisectional curvature. Using this Liouville theorem
 we prove a splitting theorem for such manifolds as well as  a gap theorem in
terms of the curvature decay of such a manifold.

In [N],
the first author raised  the following question:

\medskip

{\it On a complete noncompact K\"ahler manifold with nonnegative Ricci
curvature,
is  a plurisubharmonic function of sub-logarithmic growth  a constant?}

\medskip

\noindent It is well-known that for the complex Euclidean space $\C^m$, the
answer is positive.  An affirmative answer to the above question is also a
natural
analogue, for plurisubharmonic functions, of Yau's Liouville theorem [Y] for
positive harmonic functions on Riemannian manifolds with nonnegative Ricci
curvature.  In this paper we shall first prove the following result as a
supporting evidence of the positive solution to the above mentioned
question.

\proclaim {Theorem 1}
Let $M$ be a complete K\"ahler manifold with nonnegative holomorphic
bisectional curvature.  Let $u$ be a
plurisubharmonic function on $M$ satisfying $\D u\le \exp(C(r^2(x)$ $+1))$
 for some $C>0$. Suppose
$$
\limsup_{x\to\infty}\frac{u(x)}{\log r(x)}=0,
$$
then $u$ must be a constant.
\endproclaim

Even the result holds only for the manifolds with nonnegative bisectional
curvature it has interesting applications in studying the geometry
of such  complete K\"ahler manifolds. The first application is the
following splitting result.

\proclaim {Theorem 2}
Let $M^m$ be a complete noncompact K\"ahler manifold with bounded nonnegative holomorphic bisectional
curvature. Suppose $f$ is a nonconstant harmonic function on $M$ such that
$$
\limsup_{x\to\infty}\frac{|f(x)|}{r^{1+\e}(x)}=0, \tag 0.1
$$
for any $\e>0$, where $r(x)$ is the distance of $x$ from a fixed point.
Then  $f$ must be of linear growth and  $M$ can be splitted
  isometrically as $\widetilde{M}\times\R$.
 Moreover the universal cover $\overline {M}$ of $M$  can be splitted
  isometrically and holomorphically as $\widetilde{M'}\times\C$,
where $\widetilde M'$ is a complete K\"ahler
  manifold with nonnegative holomorphic bisectional curvature. Suppose that
there exists a holomorphic function $f$ on $M$ of growth satisfying (0.1).
Then $M$ itself splits as $\widetilde{M}\times\C$.
\endproclaim

A  consequence of Theorem 2 is the following corollary.

\proclaim{ Corollary}
Let $M^m$ be a complete noncompact K\"ahler manifold with
nonnegative holomorphic bisectional curvature whose Ricci curvature
is positive at some point. Then every  harmonic function defined on
$M$ satisfying (0.1)  must be constant.
\endproclaim

In [L1],  Li  proved that if $M^m$ is a complete noncompact K\"ahler
manifold with complex dimension $m$  with  nonnegative Ricci curvature and
$M$ supports    $n+1$  linearly independent harmonic functions of
linear growth over $\R$, then $M$ is   holomorphically isometric
to $\C^m$. Here $n=2m$ is the real dimension of $M$.
This  result was later
generalized to the real case by Cheeger-Colding-Minicozzi in [C-C-M].
They proved that Li's result   is still true for Riemannian manifolds
with nonnegative Ricci curvature. In this case, the  conclusion is that the
manifold is isometric to the Euclidean space.  In fact, they  proved that
if $M^n$ is a complete noncompact Riemannian manifold with nonnegative
Ricci curvature which supports a  non-constant harmonic function of
linear growth, then the tangent cone at infinity splits a factor of $\R$.
Theorem 2  shows that on a complete K\"ahler manifold with bounded
nonnegative holomorphic bisectional
curvature, the existence of a nonconstant
linear growth harmonic function would split the manifold itself. One can
also think this as a function-theoretic version of Cheeger-Gromoll's
splitting theorem.

On the other hand, in [Y] and [C-Y],  Cheng and Yau proved that on a
complete noncompact Riemannian manifold with nonnegative Ricci curvature,
then any sublinear growth harmonic function must be constant.
That is to say, if the growth rate of a harmonic function is `close'
to that of constant functions, then the harmonic function must be constant.
It is an interesting question to locate the `next gap'. Namely, what is the
minimum growth rate beyond the linear growth.
An easy consequence of the theorem is that if a harmonic function
 $f(x)$ is of $O\left(r(x)\left(\log r(x)\right)^a\right)$, for some $a>0$,
then  it is of linear growth.
  On the other hand, for any $\delta>0$,
the `round off' cones with metrics
$dr^2+r^2\ ds^2_{S^1(1+\delta)}$, where $S^1(\frac{1}{\sqrt{1+\delta}})$
is the circle with radius $\frac{1}{\sqrt{1+\delta}}$, support harmonic functions of growth $r^{1+\delta}(x)$. Therefore, Theorem 1 provides
the best `next gap', at least for the K\"ahler manifolds
with bounded nonnegative
bisectional curvature.
Whether the similar `gap' exists for the manifolds with nonnegative
Ricci curvature
and whether the same splitting result remains true for the Riemannian
manifolds with nonnegative sectional curvature remain to be interesting open
questions.

Theorem 1 also shows its potential in the study of the structure of complete
K\"ahler manifolds with nonnegative curvature by proving the following gap
theorem.

\proclaim{Theorem 3}
Let $M$ be a complete noncompact K\"ahler manifold with
nonnegative bisectional curvature.
Assume that  ${\Cal R}(x)\le C(r^2(x)+1)$ for some $C>0$ and
$$\int_0^rs\lf(\aint_{B_o(s)}{\Cal R}(y)\, dy\ri)ds=o(\log r)
$$
 where ${\Cal R}(x)$ is the scalar curvature
function, $\aint_{B_o(r)}{\Cal R}(y)\, dy$ is the average of $\Cal R$ over $B_o(r)$ . Then $M$ must be flat. In particular, the universal cover of $M$ must be
isometric to $\C^m$.
\endproclaim

Theorem 3 is the best gap type theorem proved so far for
the K\"ahler manifolds with nonnegative bisectional curvature.
 The first result of this
sort was proved by Mok-Siu-Yau in [M-S-Y] through
solving the Poincar\'e-Lelong
equation. It was later generalized in [N] for non-parabolic manifolds.
The best
 result to date is in
[C-Z], where Chen and
Zhu used   W.-X. Shi's argument in [Sh3] and proved  that:

\medskip

{\it Let $M^m$ be a complete noncompact K\"ahler manifold with  nonnegative
bisectional curvature. Assume that   $\Cal R(x)$ is bounded and
$\aint_{B_x(r)}{\Cal R}(y)\, dy
\le k(r) $  for all $x$ where    $k(r)$ is a nonincreasing function satisfying $k(r)=o(r^{-2})$.
Then $M$ must be flat. }

\medskip

\noindent
In the proof of [C-Z], the long time existence on
the K\"ahler-Ricci flow in  [Sh3] was used together
with the volume element estimate in [Sh3] and a Li-Yau-Hamilton
inequality on K\"ahler-Ricci flow of H.-D. Cao [Co1].
Along this line, the authors of the current
paper improved the above result slightly in [N-T2], after
 their simpler derivation of W.-X. Shi's volume element estimate.

We prove Theorem 3  using a simpler and much more direct method.
Note that in Theorem 3
we do not require the uniform decay as in the above mentioned
result in [C-Z] (See also [N-T2]).
Another advantage of this method is that we do not necessarily require
the boundedness of the
curvature tensor, which has to be assumed due to the current status
 of the existence theory on the K\"ahler-Ricci flow over complete
 noncompact manifolds.
The connection between the Liouville theorem and the gap theorem is provided
by the solution to the Poincar\'e-Lelong equation, especially the one
constructed in Theorem 5.1  of [N-S-T1]. This connection was
 also illustrated in earlier papers [N-T1-2].  We should also remark that
Theorem 3 is not true if one only assume that the manifold has nonnegative
Ricci curvature. In fact, there are many examples of K\"ahler manifolds
with curvature decay satisfying Theorem 3 and with maximum volume growth.
But they are not flat. (To our knowledge that all such examples are Ricci
flat. Whether this is generally true or not remains an interesting question.)

Finally,  we summarize some  earlier work on Liouville properties of the
plurisubharmonic functions and explain the methods we use to
 put our theorem into the right
perspective. The first attempt to the question raised at the beginning on the
plurisubharmonic functions was made in [N], where the first author proved
that if $M$ is quasi-projective and $u$ is bounded then $u$ is a constant.
In general case, it was also proved there that
$u$  satisfies a homogenous Monge-Amper\'e equation. This
fact turns out to be very useful in our proof here.
In [N-T1], using the
K\"ahler-Ricci flow and a linear trace
Li-Yau-Hamilton inequality (which is also
called Harnack inequality in [H2], [Co1] and [C-H])
established there, the authors  answered the question raised at the beginning
affirmatively
under the assumptions that $M$ has {\it bounded} bisectional
curvature,
the average of the scalar curvature of $M$ has quadratic decay
and the Laplacian
of the plurisubharmonic function grows at most exponentially (cf.
 [Theorem 3.2, N-T1] for a more precise statement).

The proof here
is complete different, simpler and is based on the fact
that on a complete K\"ahler
manifold with nonnegative Ricci curvature, a plurisubharmonic function of
sub-logarithmic growth must satisfies a homogeneous Monge-Amper\'e equation
(see Lemma 2.3 below).
Namely, at each point, at least one of the eigenvalue of the complex
Hessian of the plurisubharmonic function is zero. Hence a
natural way to prove the theorem is to use the induction. However, the
foliation defined by the Monge-Amper\'e equation might be singular. We
overcome this difficulty by deforming the plurisubharmonic function through
heat equation. It turns out that, under the condition that the manifold has
nonnegative holomorphic bisectional curvature, the deformed function is
still plurisubharmonic for  $t>0$ and satisfies the homogeneous
Monge-Amper\'e
equation.  Then the manifold, or
its universal cover if not simply-connected, can be splitted with a factor
whose tangent space corresponds to the kernel of the complex Hessian of the
function at
each point at some time $t>0$. (Namely, the foliation at $t>0$
becomes a  product.) Therefore, we can
indeed use the induction to conclude
the result.

We should point out that
the upper-bound assumption on the Laplacian is believed not necessary.
However, due to the lack of direct method to the problem and the
 heat  equation method in our proof, this assumption is necessary
to obtain a  maximum principle which
holds for tensors satisfying a heat equation.
The classical uniqueness for the solution to the heat equation on Euclidean
space requires a similar necessary assumption on the solution.

\medskip

{\bf Acknowledgment.}
The authors would like to thank Huai-Dong Cao, Bing-Long Chen, Bennet Chow,
Jiaping Wang and Fangyang Zheng for helpful discussions.
In particular, the discussion with Fangyang Zheng leads to  an
improvement of Theorem 2 and its consequence.
We also thank Peter Li, Richard Schoen
 for their   interests.

\input amstex
\documentstyle{amsppt}
\magnification=1200
\hsize=13.8cm
\catcode`\@=11
\def\NoLogo{\let\logo@\empty}
\catcode`\@=\active
\NoLogo

\def\heat{\lf(\frac{\p}{\p t}-\Delta\ri)}
\def \b {\beta}
\def\i{\sqrt{-1}}

\def\lf{\left}
\def\ri{\right}
\def\bbar{\bar \beta}
\def\a{\alpha}
\def\ol{\overline}
\def\g{\gamma}
\def\e{\epsilon}
\def\p{\partial}
\def\delbar{{\bar\delta}}

\def\dbar{\bar\partial}

\def\C{\Bbb C}
\def\R{\Bbb R}

\def\vp{\varphi}

\def\dbar{\bar\partial}
\def\ba{{\bar\alpha}}
\def\bb{{\bar\beta}}
\def\abb{{\alpha\bar\beta}}

\def\i{\sqrt {-1}}

\def \D {\Delta}
\def\aint{\frac{\ \ }{\ \ }{\hskip -0.4cm}\int}
\magnification=1200

\subheading{\S1 Estimates on solutions to the heat equation}

In this section, we derive   some basic estimates on  the solution  of the heat
equation with plurisubharmonic initial data
on a complete noncompact K\"ahler manifold  with nonnegative bisectional curvature.
First, we need the following:
\proclaim{Lemma 1.1} Let $M^n$ be a complete noncompact Riemannian manifold
with nonnegative Ricci curvature. Let $u_0$ be a smooth function on $M$
such that
$$ \exp\lf(a(r^2(x)+1)\ri)\le  u_0(x) \le  \exp\lf(b(r^2(x)+1)\ri)
$$
for some constants $b>a>0$, where $r(x)$ is the distance of $x$ from a fixed
point $o\in M$. Then there exists a $T>0$ depending only on $b$ such that
the Cauchy problem
$$
\heat u=0,\qquad u(x,0)=u_0(x)\tag1.1
$$
has a solution $u$ on $M \times [0,T]$. Moreover, there exist constants
$C_1,\ C_2>0$ such that
$$ C_1\exp\lf(\frac a4 r^2(x)\ri)\le u(x,t)\le  C_2\exp\lf( 3br^2(x)\ri) \tag 1.2
$$
on $M\times[0,T]$.
\endproclaim
\demo{Proof}  To prove the existence of the  solution, it is sufficient to show that for some $T>0$,
$$
u(x,t)=\int_M H(x,y,t) u_0(y)dy
$$
is uniformly bounded on $B_o(r)\times[0,T]$ for all $r$. Let $x\in M$ and let   $r(x)=2R$, using the  upper bound of the heat kernel of Li-Yau \cite{L-Y}, we have that
$$
\split
\int_M H(x,&y,t)\,  u_0(y)dy\\
&=\int_{B_x(R)}H(x,y,t)u_0(y)dy+\int_{M\setminus B_x(R)}H(x,y,t)u_0(y)dy\\
&\le  \exp\lf(b(9R^2 +1)\ri)+ C_3 \int_{M\setminus B_x(R)} V_x^{-1}(\sqrt t)  \exp\lf(-\frac{r^2(x,y)}{5t}+9b r^2(x,y)\ri)dy  \endsplit
$$
for some constant $C_3$ depending only on $n$ and $b$, where we have used the fact that $\int_M H(x,y,t)dy=1$ and the fact that $r(y)\le 3r(x,y)$ outside $B_x(R)$. Here $r(x,y)$ is the distance between $x$ and $y$.
  If we choose $T =\min\{1,\frac{1}{100b}\}$, then for $0<t\le T$,
 we have that
$$
\split
\int_M H(x,y,t) u_0(y)dy
&\le \exp\lf( b(9R^2 +1)\ri)+\frac{C_4R^n}{V_o(1)}\int_{R}^\infty   t^{-\frac n2}\exp\lf(-\frac{ r^2}{10t}\ri)r^{n-1}dr\\
&\le \exp \lf(b(9R^2+1)\ri)+C_5R^{n}\\
& \le C_6  \exp \lf(b(9R^2+1)\ri)
\endsplit
$$
for some constants $C_4- C_6$ depending only on $n$, $b$ and $V_o(1)$.
Here we have used the
  volume comparison so that $V_x(1)\ge  \frac{1}{(1+2R)^n}V_o(1)$ and
$A(\p B_x(r))\le C(n)r^{n-1}$. From this it is easy to see that (1.1) has a solution $u(x,t)$ so that the second inequality of (1.2) is true. To prove the lower bound of $u(x,t)$,   if $2R=r(x)>2$, then for $0<t\le T$,
$$
\split
u(x,t)&\ge \int_{B_x(R)}H(x,y,t)u_0(y)dy\\
&\ge C_7\exp \lf(aR^2\ri)\int_{B_x(\sqrt t)}V^{-1}_x(\sqrt t) \exp\lf(-\frac{r^2(x,y)}{3t}\ri)dy\\
&\ge C_8\exp\lf(aR^2\ri)
\endsplit
$$
for some positive constants $C_7$ and $C_8$ depending only on $n$ and $a$, where we have used the lower estimate of the heat kernel in \cite{L-Y,
page 182}. From this it is easy to see the first inequality in (1.2) is also true.
\enddemo
Similarly, one can prove that:
\proclaim{Lemma 1.2}
 Let $M^n$ be a complete noncompact manifold with nonnegative Ricci
curvature and let $u_0(x)$ be a smooth function on $M$ such that
$|u_0(x)|\le \exp\lf(b(r^2(x)+1)\ri)$, then the Cauchy problem (1.1) has a solution   $u(x,t)$ on $M\times [0,T]$ for some $T>0$, such that
$|u(x,t)|\le C\exp\lf(3b r^2(x)\ri)$ for all $(x,t)$. If in addition,
$$
\limsup_{x\to\infty}\frac{u_0(x)}{\log r(x)}=0\tag 1.3
$$
then for any $T\ge t>0$,
$$
\limsup_{x\to\infty}\frac{u(x,t)}{\log r(x)}=0.\tag1.4
$$
\endproclaim
In the next lemma, we give  sufficient conditions  for a function $u(x)$ to satisfy $|u(x)|\le \exp\lf(C(r^2(x)+1)\ri)$.

\proclaim{Lemma 1.3} Let $M^n$ be a complete noncompact manifold with nonnegative Ricci curvature. Let $u$ be a smooth function. Suppose  there exists a constant $b>0$ such that
\roster
\item"{(i)}" $u(x)\le \exp\lf(b(r^2(x)+1)\ri)$;
 and
\item"{(ii)}"  $0\le \Delta u(x)\le \exp\lf(b(r^2(x)+1)\ri)$.
\endroster
 Then
$$
u(x)\ge - C\exp\lf(5b(r^2(x)+1)\ri)\tag1.5
$$
for some constant $C>0$ for all $x$.
\endproclaim
\demo{Proof} Consider $\tilde M= M\times \Bbb R^4$, then $\tilde M$ has positive Green's function $G(x,y)$ which satisfies:
$$
G(\tilde x,\tilde y)\le \frac{C\tilde r^2(\tilde x,\tilde y)}{V_x\lf(\tilde r(\tilde x,\tilde y)\ri)}\tag1.6
$$
for some constant $C$ depending only on $n$, see \cite{Sh3, p.162}. If we define $\tilde u(\tilde x)=u(x)$, where $\tilde x=(x,x')\in \tilde M = M\times \Bbb R^4$, then $\tilde u$ also satisfies conditions similar to (i) and (ii). Suppose we can prove that (1.5) is true for $\tilde u(\tilde x)$ on $\tilde M$, then it is easy to see that $u$ also satisfies (1.5). Hence we may assume that $M$ has a positive Green's function which satisfies condition similar to (1.6) on $M$.

Now for any $R>0$, let $f(x)=\Delta u(x)$ and let
$$
v_R(x)=-\int_{B_o(R)}G_R(x,y)f(y)dy.
$$
Then for $x\in B_o(\frac R2)$,
$$
\split
|v_R(x)| &\le \exp\lf(b(R^2 +1)\ri)\int_{B_o(R)}G(x,y)dy\\
&\le \exp\lf(b(R^2 +1)\ri)\int_{B_x(2R)}G(x,y)dy\\
&\le C_1\exp\lf(b(R^2 +1)\ri)\int_0^{2R}\frac{r^2A_x(r)}{V_x(r)}dr\\
&\le C_2R^2\exp\lf(b(R^2 +1)\ri)
\endsplit\tag1.7
$$
for some constants $C_1-C_2$ independent of $x$ and $R$. Here we have used condition (ii), (1.6) and volume comparison. Since $v_R$ satisfies $\Delta v_R=f$ on $B_o(R)$ with zero boundary data, we conclude that $v_R+\exp\lf(b(R^2 +1)\ri)\ge u$ in $B_o(R)$ by condition (i) and the maximum principle. By the Harnack inequality derived from the gradient estimate for positive harmonic functions of Cheng-Yau \cite{C-Y}, we have
$$
\sup_{B_o(\frac R2)} \lf(v_R+\exp\lf(b(R^2 +1)\ri)-u\ri)\le C_3\lf(v_R(o)+\exp\lf(b(R^2 +1)\ri)- u(o)\ri)
$$
for some constant $C_3$ depending only on $n$. Hence for $x\in B_o(\frac R2)$,
$$
\split
-u(x)&\le -v_R(x)+C_3\lf( \exp\lf(b(R^2 +1)\ri)- u(o)\ri)\\
&\le C_2R^2\exp\lf(b(R^2 +1)\ri)+C_3\lf( \exp\lf(b(R^2 +1)\ri)- u(o)\ri)
\endsplit
$$
where we have used (1.7) and the fact that $v_R\le 0$. From this, it is easy to see that (1.5) is true.
\enddemo

Let  $M^m$ be a complete K\"ahler manifold and
let $u(x,t)$ be a solution to the heat equation on $M$. Namely,
$$
\left(\frac{\p}{\p t}-\D\right) u(x,t)=0.
$$
By the computation in   Lemma 2.1 in \cite{N-T1}, we have:
\proclaim{Lemma 1.4}
Let $u(x,t)$ be a solution to the heat equation. Then the complex Hessian $u_{\abb}(x,t)$ of $u(x,t)$ satisfies the
complex Lichnerowicz equation:
$$
\left(\frac{\p}{\p t}-\D \right)u_{\g\delbar}=
R_{\beta \bar{\a}\g\delbar}u_{\a\bbar}-
\frac{1}{2}\left(R_{\g\bar{p}}u_{p\delbar}+
R_{p\delbar}u_{\g\bar{p}}\right).\tag1.8
$$ Here $R_{\beta \bar{\a}\g\delbar}$ and $R_{\g\bar{p}}$ are the
  curvature tensor and the Ricci tensor of $M$.
\endproclaim

\proclaim{Lemma 1.5}
Let $M^m$ be a complete noncompact K\"ahler manifold with nonnegative
bisectional curvature. Let $u(x,t)$ be a solution to the heat equation. Then
$$
\left(\frac{\p}{\p t}-\D \right)\|u_{\abb}\|^2 \le -\|\nabla_\g u_\abb\|^2
-\|\nabla_{\bar\g}u_\abb\|^2. \tag 1.9
$$
\endproclaim
\demo{Proof} Choose a normal coordinate. Using Lemma 1.4, the
 direct calculation shows that:
$$
\split
\D \|u_{\abb}\|^2& =\|\nabla_\g u_\abb\|^2
+\|\nabla_{\bar\g}u_\abb\|^2+u_{\abb}\overline{u_{\g\bbar}}R_{\a\bar{\g}s\bar{s}}
+u_{\abb}\overline{u_{\a\bar{\g}}}R_{\g\bbar s \bar{s}}\\
& \ \ -\overline{u_{\abb}}u_{s\bar{t}}R_{\abb s \bar{t}}-u_{\abb}
\overline{u_{s\bar{t}}R_{\abb s \bar{t}}}+u_{\abb}\overline{(u_t)_{\abb}}
+\overline{u_{\abb}}(u_t)_{\abb}.
\endsplit
$$
Therefore,
$$
\left(\frac{\p}{\p t}-\D \right)\|u_{\abb}\|^2=-\|\nabla_\g u_\abb\|^2
-\|\nabla_{\bar\g}u_\abb\|^2-\sum_{\a\b}R_{\a\bar{\a}\b\bar{\b}}(\lambda_{\a}-
\lambda_{\b})^2,
$$
where $\lambda_{\a}$ are eigenvalues of $u_{\abb}$. Here we calculate under a normal
coordinate around a fixed point such that $u_{\abb}$ is diagonalized. By the
nonnegativity of the bisectional curvature the lemma follows.
\enddemo

\proclaim{Lemma 1.6} Let $M^m$ be a complete noncompact K\"ahler manifold with nonnegative bisectional curvature.
Let $u_0$ be a plurisubharmonic function on $M$ such that
$$
 u_0(x) \le \exp\lf(b(r^2(x)+1)\ri),\tag1.10
$$
and
$$
0\le \Delta u_0(x)\le \exp\lf(b(r^2(x)+1)\ri)\tag1.11
$$
for some constant $b>0$. Then there exists $T>0$ such that the Cauchy problem (1.1) has a solution $u(x,t)$ on $M\times[0,T]$
such that for some constant $b^*>0$
(which might depend on $T$),
$$
||u_\abb||(x,t)\le \exp\lf(b^*(r^2(x)+1)\ri) \tag 1.12
$$
for all $(x,t)$.
\endproclaim
\demo{Proof} By Lemmas 1.2 and 1.3, we conclude that for some $T>0$, (1.1) has a solution $u(x,t)$ in $M\times[0,T]$ such that
$$
|u(x,t)|\le C_1\exp\lf(15br^2(x)\ri).\tag1.13
$$
It remains to prove (1.12). By (1.11) and (1.13), one can easily prove that
$$
\int_{B_o(r)}|\nabla u_0|^2 dx\le \exp\lf(b_1(r^2+1)\ri)\tag1.14
$$
for some constant $b_1$ for all $r$.
Since
$$
\heat u^2=-2|\nabla u|^2,
$$
we can multiply the equation by $\vp^2(x)$ and integrate by parts, where $\vp(x)$ is a smooth function such that $0\le \vp \le 1$,  $\vp=1$ on
$B_o(R)$, $\vp=0$ outside $B_o(2R)$ and $|\nabla\vp|\le C/R$ for some constant $C$ independent of $R$. We have
$$
\split
2\int_0^T\int_M \vp^2|\nabla u|^2 dxdt&=-\int_0^T\int_M \vp^2\heat u^2\\
&\le \int_M \vp^2 u_0^2(x) dx+4\int_0^T\int_M  \vp u|\nabla \vp|\,|\nabla u|dxdt\\
&\le \int_M \vp^2 u_0^2(x) dx+4\int_0^T\int_M |\nabla \vp|^2u^2dxdt+\int_0^T\int_M \vp^2|\nabla u|^2 dxdt.
\endsplit
$$
Hence by  (1.13), we have
$$
\int_0^T\int_M \exp\lf(-b_2(r^2(x)+1)\ri)|\nabla u|^2(x)dxdt<\infty \tag1.15
$$
for some $b_2>0$. As in \cite{N-T2, Lemma 1.1}, using the fact that the Ricci curvature is nonnegative, we have
$$
\heat |\nabla u|^2\le  -  ||u_\abb||^2-||u_{\a\b}||^2.
$$
Multiplying  this by $\vp^2$ and integrating by parts, using (1.14)
and (1.15),
one can repeat the above argument and show that

$$
\int_0^T\int_M \exp\lf(-b_3(r^2(x)+1)\ri)||u_\abb||^2(x)dxdt<\infty \tag1.16
$$
for some constant $b_3>0$.

To conclude the proof of (1.12), let $\Phi=||u_\abb||^2$, by Lemma 1.5  we have
$$
\heat \Phi(x,t)\le 0.\tag1.17
$$
  By the mean value inequality \cite{Theorem 1.2,  L-T},
we have, for $  \frac 1 4 r^2\ge T$,
$$
\exp(-C_1r^2t)(1+\Phi)(x,t) \le C_2\bigg[ \frac1{r^{2m+2}}\int_0^T\int_{B_x(\frac r4)}\exp(-C_1r^2s) \Phi (y,s)dyds +\sup_{B_x(\frac r4)}\Phi(y,0)\bigg].
$$
Combining this with (1.11) and (1.16),  we can conclude that (1.12) is
true.
\enddemo

\input amstex 
\documentstyle{amsppt}
\magnification=1200
\hsize=13.8cm
\catcode`\@=11
\def\NoLogo{\let\logo@\empty}
\catcode`\@=\active
\NoLogo

\def\heat{\lf(\frac{\p}{\p t}-\Delta\ri)}
\def \b {\beta}
\def\i{\sqrt{-1}}

\def\lf{\left}
\def\ri{\right}
\def\bbar{\bar \beta}
\def\a{\alpha}
\def\ol{\overline}
\def\g{\gamma}
\def\e{\epsilon}
\def\p{\partial}
\def\delbar{{\bar\delta}}

\def\dbar{\bar\partial}

\def\C{\Bbb C}
\def\R{\Bbb R}

\def\vp{\varphi}

\def\dbar{\bar\partial}
\def\ba{{\bar\alpha}}
\def\bb{{\bar\beta}}
\def\abb{{\alpha\bar\beta}}

\def\i{\sqrt {-1}}

\def \D {\Delta}
\def\aint{\frac{\ \ }{\ \ }{\hskip -0.4cm}\int}
\documentstyle{amsppt}
\magnification=1200
\hsize=13.8cm
\subheading{\S2 Proof of the Liouville theorem}

The proof of Theorem 1 is
based on the following three lemmas. Let $M^m$ and $u_0$ as in Lemma 1.6, and let $u(x,t)$ be the solution of (1.1) on $M\times[0,T]$  for some $T>0$ constructed in the lemma. In the following, the eigenvalues of a Hermitian form are arranged in ascending order. Hence the first eigenvalue is the smallest one.
\proclaim{Lemma 2.1} With the above notations and assumptions, we have the following:
\roster
\item"{(a)}" $u_\abb(x,t)\ge 0$ for all $(x,t)$.
\item"{(b)}" For any $T>t'\ge0$, suppose there is a point $x'$ in $M^m$ such that the first  $k$   eigenvalues $\lambda_1,\dots,\lambda_k$ of $u_\abb(x',t')$ satisfy $\lambda_1+\dots+\lambda_k>0$, then for all $t>t'$ and for all $x\in M$, the sum of the first $k$  eigenvalues of $u_\abb(x,t)$ is also positive.
\endroster
\endproclaim
 \demo{Proof}   The proofs of (a) and (b) are similar. Let us first prove (b)  for the case that $t'=0$.   By Lemma 1.6, there is a constant $b_1>0$ such that on $M\times[0,T]$
$$
||u_\abb||(x,t)\le \exp\lf(b_1(r^2(x)+1)\ri).\tag2.1
$$
It is easy to see that there exists a smooth function $h_0(x)$ such that
$$
\exp\lf(8b_1(r^2(x)+1)\ri)\le h_0(x)\le \exp\lf(b_2 (r^2(x)+1)\ri)
$$
for some $b_2>8b_1$. By Lemma 1.1, we can find a solution $h(x,t)$ of the heat equation on $M\times[0,T_1]$ where $T_1=\min\{1,\frac{1}{100b_2}\}$ such that
$$
h(x,t)\ge C_1\exp (2b_1r^2(x)) \tag2.2
$$
for some $C_1>0$. Let $\phi(x,t)=e^t h(x,t)$, then
$$
\heat \phi=\phi.\tag2.3
$$
Assume at $t=0$, there is $x_0\in M$, such that the sum of the first $k$   eigenvalues of $u_\abb(x_0,0)$ is positive. Then we can find a smooth nonnegative function $f$ with $f(x_0)>0$ and with support in a neighborhood of $x_0$, such that the sum of the first $k$   eigenvalues of
$$
u_\abb-fg_\abb
$$
is still nonnegative at $t=0$, where we have used the fact that $u_\abb\ge 0$ at $t=0$. As in \cite{B}, solve the scalar heat equation
$$
\heat f=-f\tag 2.4
$$
such that $f(x,0)=f(x)$. The solution can simply be written as
$e^{-t}\cdot\int_M H(x,y,t)f(y)\, dy$.
Then by the strong maximum principle, $f>0$ for $t>0$ and $f$ is bounded.

Let $\e>0$, define $\psi=-f+\e\phi$, and let $\eta_\abb=u_\abb+\psi g_\abb$, where $g_\abb$ is the metric tensor of $M$. Then at $t=0$, at each point the sum of the first $k$  eigenvalues of $\eta $ is positive. We want to prove that for any $T_1\ge t>0$ and any $x\in M$, the sum of the first $k$   eigenvalues of $\eta$ is positive. Otherwise, because of (2.1), (2.2), the definition of $\phi$ and the fact that $f$ is bounded,  it is easy to see that there is a first $T_1\ge t_1>0$ and a point $x_1\in M$, such that the sum of the first $k$ eigenvalues of $\eta$ at $x_1$ at time $t_1$ is zero.

Let us fix the notations. Suppose $v_1,\dots,v_m$ are unit
eigenvectors of $\eta$ at $x_1$ at time $t_1$, with eigenvalues
$\lambda_1\le \lambda_2\le\dots\le \lambda_m$.   We may choose   normal
coordinates at $x_1$   such that $v_j=\frac{\p}{\p z^j}$ at $x_1$.
In particular, if we write $v_j=v_j^\a\frac{\p}{\p z^\a}$, we have $v_j^\a=\delta_{\a j}$ at $x_1$.  Note  that the sum of the first $k$ eigenvalues of a Hermitian form is   the infimum of  the traces of the form  restricted to
  $k$-dimensional subspaces. Therefore $\sum_{\a, \beta =1}^{k}
\lf(g^{\abb}\eta_{\abb}\ri)\ge 0$ for all $(x,t)$ with $t\le t_1$ and
equals to zero at $(x_1, t_1)$.

Hence at $(x_1,t_1)$, we have
$$
0\ge \heat\lf(\sum_{ \a,\beta =1 }^k  \eta_\abb g^{\abb}\ri).\tag 2.5
$$
From now on repeated indices mean summation from $1$ to $m$ if there is no specification.
Now
$$
\frac{\p}{\p t}\lf(\sum_{\a,\beta =1}^k  \eta_\abb g^{\abb}\ri)
=\sum_{\a,\beta=1}^k\lf(\frac{\p}{\p t}\eta_\abb \ri)g^{\abb}.
\tag2.6
$$
Also at $(x_1,t_1)$, we have
$$\Delta \lf(\sum_{\a,\beta=1}^k  \eta_\abb g^{\abb}\ri)
= \sum_{\a,\beta=1}^k \lf(\Delta\eta_\abb\ri)    g^{\abb}.
$$
Combining this with (2.5), (2.6) and (1.8) in Lemma 1.4,   at $(x_1,t_1)$ we have,
$$
\split
0&\ge \sum_{\a,\beta=1}^k \bigg[R_{\delta\bar\g\a\bb}\lf(u_{\g\delbar}+\psi g_{\g\delbar}\ri)-\frac12R_{\a\bar p}\lf(u_{p\bb}+\psi g_{p\bb}\ri) -\frac12R_{p\bb }\lf(u_{\a\bar p }+\psi g_{\a\bar p }\ri)\bigg]g^{\abb}\\
&+\sum_{\a,\beta=1}^k\bigg(\lf[\heat\psi\ri]g_\abb -R_{\delta\bar\g\a\bb}\psi g_{\g\delbar}+\frac12\psi R_{\a\bar p}g_{p\bb} +\frac12\psi R_{\a\bar p }g_{p\bb}\bigg)g^{\abb} .\endsplit\tag2.7
$$
Since at $(x_1,t_1)$, $\eta$ has eigenvectors $v_p=\frac{\p}{\p z^p}$, for $1\le p\le m$, with eigenvalue $\lambda_p$
$$
\split
\sum_{\a,\beta=1}^k
\bigg[R_{\delta\bar\g\a\bb}&\lf(u_{\g\delbar}+\psi g_{\g\delbar}\ri)-\frac12R_{\a\bar p}\lf(u_{p\bb}+\psi g_{p\bb}\ri) -\frac12R_{p\bb }\lf(u_{\a\bar p }+\psi g_{\a\bar p }\ri)\bigg]g^{\abb}\\
&=\sum_{j=1}^k\sum_{\gamma=1}^m R_{\gamma\bar\g j\bar j}\lambda_\g-\sum_{j=1}^k R_{j\bar j}\lambda_j\\
&=\sum_{j=1}^k\sum_{\gamma=1}^m R_{\gamma\bar\g j\bar j}\lambda_\g-\sum_{j=1}^k \sum_{\gamma=1}^mR_{\g\bar\g j\bar j}\lambda_j\\
&=\sum_{j=1}^k\sum_{\g=k+1}^m\lambda_\g R_{\g\bar \g j\bar j}-\sum_{j=1}^k\sum_{\g=k+1}^mR_{\g\bar \g j\bar j}\lambda_j\\
&=\sum_{j=1}^k\sum_{\g=k+1}^m R_{\g\bar \g j\bar j}(\lambda_\g-\lambda_j)\\
&\ge 0
\endsplit\tag2.8
$$
where we have used that fact that $M$ has nonnegative bisectional curvature, and $\lambda_\g\ge \lambda_j$ for $\g\ge j$. Also by (2.3) and (2.4)
$$
\lf[\heat\psi\ri]=k(f+\e\phi)>0.
\tag2.9
$$
Moreover
$$
 \lf(-R_{\delta\bar\g\a\bb}\psi g_{\g\delbar}+\frac12\psi R_{\a\bar p}g_{p\bb} +\frac12\psi R_{\a\bar p }g_{p\bb}\ri)g^{\abb}=0.\tag2.10
$$
From (2.7)--(2.10), we have a contradiction. Hence the sum of the first $k$ eigenvalues of $\eta$ is nonnegative for all $(x,t)\in M\times(0,T_1]$. Let $\e\to0$, we conclude that   the sum of the first $k$   eigenvalues of   $u_\abb(x,t)-f(x,t)g_\abb(x,t)$ is nonnegative on $M\times[0,T_1]$. Since $f$ is positive for $t>0$, the  sum of the first $k$  eigenvalues of   $u_\abb(x,t)$ must be positive for $0<t\le T_1$.

Take $f\equiv0$ in the above, one can prove similarly that $u_\abb\ge0$ on $M\times[0,T_1]$. One can then apply the same arguments as above to prove that (b) is true on $[0,T_1]$. The results then follow by iteration, because $T_1>0$ is a fixed number.

\enddemo

\proclaim{Lemma 2.2} Let $M$,  $u_0(x)$ and $u(x,t)$ be
as in Lemma 2.1.  Let
$$
\Cal K(x,t)=\{v\in T^{1,0}_x(M)|\  u_\abb(x,t) v^\a=0,\text{\rm \ for all $\b$}\}
$$
be the null space of $u_\abb(x,t)$. Then there exists $0<T_1<T$ such that  for any $0<t<T_1$, $\Cal K(x,t)$ is   distribution on $M$. Moreover the distribution is
invariant under parallel translation.
\endproclaim
\demo{Proof}  By Lemma 2.1, it is easy to see that there exists $0<T_1<T$ such that $\dim \Cal K(x,t)$ is constant on $M\times(0,T_1)$. Hence for each $0<t<T_1$, $\Cal K(x,t)$ is a smooth distribution on $M$. It remains to prove that the distribution is parallel for fixed $t$. We can proceed as in \cite{H1, Lemma 8.2}.

Fix $0<t_0<T_1$, let $x_0\in M$ and let $w_0\in \Cal K(x_0,t_0)$. Let $\gamma(\tau)$ be a smooth curve from $x_0$ and let $w(\tau)$ be the vector field obtained by parallel translation along $\gamma$. We want to prove that $w(\tau)$ is also in the null space  $\Cal K(\gamma(\tau),t_0)$ at $\gamma(\tau)$. First extend $w$ to be a vector field in a neighborhood of $\gamma(\tau)$, and then extend $w$ to be a vector field independent of time $t$. Now, projecting $w$ onto $\Cal K(x,t)$, we have a vector field $v$ such that  $v$ is in $\Cal K(x,t)$ for all $x$ in a neighborhood of $\gamma$ and for all  $t$. The following computations are performed in a neighborhood of $\gamma$.

Since
$$
u_{\abb}v^\a=0\tag2.10
$$
for all $\b$, we have
$$
\split
0&=\frac{\p}{\p t}\lf(u_{\abb}v^\a\ol{v^\b}\ri)\\
&=\lf(\frac{\p}{\p t} u_{\abb}\ri)v^\a\ol{v^\b}+u_{\abb}\frac{\p v^\a}{\p t}\ol{v^\b}+u_{\abb} v^\a \frac{\p\ol{v^\b}}{\p t}\\
&=\lf(\frac{\p}{\p t} u_{\abb}\ri)v^\a\ol{v^\b}
\endsplit\tag2.11
$$
where we have used (2.10). Choosing a unitary frame  $e_s$ at   a point $\gamma(\tau)$, we have
$$
\split
0&=\Delta\lf(u_{\abb}v^\a\ol{v^\b}\ri)\\
&=\frac12\lf(\nabla_s\nabla_{\bar s}+\nabla_{\bar s}\nabla_s\ri)\lf(u_{\abb}v^\a\ol{v^\b}\ri)\\
&=\lf(\Delta u_{\abb}\ri)v^\a\ol{v^\b}-u_{\abb}\nabla_{\bar s}v^\a \nabla_s\ol{v^\b}- u_{\abb}\nabla_{  s}v^\a \nabla_{\bar s}\ol{v^\b}
\endsplit\tag2.12
$$
where we have used (2.10) so that
$$
\lf(\nabla_{  s}u_{\abb}\ri)v^\a=- u_{\abb}\nabla_{  s}v^\a, \ \lf(\nabla_{\bar  s}u_{\abb}\ri) v^\a=- u_{\abb}\nabla_{\bar  s}v^\a
$$
and their complex conjugates.

Combining with (1.8), (2.11), (2.12), we have
$$
0=R_{t\bar s\abb}u_{s\bar t}v^\a\ol{v^\b}+2u_{\abb}\nabla_{\bar s}v^\a \nabla_s\ol{v^\b}+2 u_{\abb}\nabla_{  s}v^\a \nabla_{\bar s}\ol{v^\b}.\tag 2.13
$$
We may choose $e_s$ so that at a point    $u_{s\bar t}=a_s\delta_{st}$. Then
$$
R_{t\bar s\abb}u_{s\bar t}v^\a\ol{v^\b}=R_{s\bar s \abb }a_sv^\a \ol{v^\b}=a_sR_{s\bar s v\bar v}\ge0
$$
because $a_s\ge0$  and $M$ has nonnegative bisectional curvature. Hence (2.13) implies that $\nabla_s v$ and $\nabla_{\bar s}v$ are in the null space  $\Cal K(\gamma(\tau),t_0)$.

Since $w(\tau)$ is parallel along $\gamma(\tau)$, and $w=v+w^\perp$, where $w^\perp$ is perpendicular to $\Cal K(\g(\tau),t_0)$, we have
$$
0=\frac{D}{d\tau}w=\frac{D}{d\tau}v+\frac{D}{d\tau}w^\perp.
$$
Hence
$$
\frac{D}{d\tau}w^\perp=-\frac{D}{d\tau}v
$$
which is in $\Cal K$.

Now
$$
\frac{d}{d\tau}\langle w^\perp,w^\perp\rangle=\langle \frac{D}{d\tau}w^\perp,w^\perp\rangle+\langle w^\perp,\frac{D}{d\tau}w^\perp\rangle=0
$$
because $\frac{D}{d\tau}w^\perp$ is in $\Cal K$ and $w^\perp$ is perpendicular to $\Cal K$. At $\gamma(0)=x_0$, $w=v_0$ and so $w^\perp=0$ at $\gamma(0)$. Hence $w^\perp=0$ for all $\tau$ and so $w$ is in $\Cal K$.
\enddemo

\proclaim{Lemma 2.3} (\cite{N, Proposition 4.1})
Let $M^m$ be a complete noncompact K\"ahler manifold of complex dimension $m$,
with nonnegative Ricci curvature. Let $u(x)$ be a plurisubharmonic function on $M$ satisfying:
$$
\limsup_{x\to\infty}\frac{u(x)}{\log r(x)}=0.\tag2.14
$$
Then $(\p\dbar u)^m =0$
\endproclaim

Proposition 4.1 stated in [N] is under the assumption that $M$ is
nonparabolic. However, the proof without any changes also works for general
complete K\"ahler manifolds with nonnegative Ricci curvature.

We are ready to prove Theorem 1.
\demo{Proof of Theorem 1} Let $M$ and $u$ satisfy the conditions in Theorem 1.  Let $\tilde M$ be the universal cover of $M$, then the distance function in $\tilde M$ dominates the distance function in $M$. Hence $\tilde M$ and the lift $\tilde u$ of $u$ also satisfy the conditions in the theorem. Therefore, we may assume that $M$ is simply connected.

By Lemmas 1.6, 2.1 and 2.2, there exists $T>0$ such that the Cauchy problem (1.1) has a solution $u(x,t)$ with $u(x,0)=u(x)$. Moreover, let $0<t_0<T$ be fixed, we have $u_\abb(x,t_0)\ge0$ and the null space $\Cal K(x,t_0)$ is a parallel distribution on $M$. Using the De Rham decomposition (Cf. Theorem 8.1, page 172 of [K-N])
we know that $M=M_1\times M_2$ isometrically and holomorphically such that $u_\abb$ is zero when restricted on $M_1$ and $u_\abb$ is positive everywhere when restricted on $M_2$. By Lemma 1.2, we still have
$$
\limsup_{x\to\infty}\frac{u(x,t_0)}{\log r(x)}=0.\tag2.15
$$
Hence when restricted on $M_2$, (2.15) is still true. This contradicts Lemma 2.3 and the fact that $u_\abb$ is positive when restricted on $M_2$, unless $M=M_1$. Hence $u_\abb(x,t_0)\equiv0$ on $M$ for all $0<t_0<T$. From this
 it is easy to see that $u_\abb(x)\equiv0$. By applying
the gradient estimate of Cheng-Yau \cite{C-Y} and the fact that
$u$ satisfies (2.14),  we know that $u$ is a constant.
\enddemo

\proclaim{Remarks} 1. There is a   result by  Cao \cite{Co2}
related to the splitting phenomena in the above proof.
Cao has told the second author that using the K\"ahler-Ricci flow he has
proved that if $M$ is a complete noncompact simply connected K\"ahler manifold
of bounded and nonnegative holomorphic bisectional curvature, then
$M$ splits holomorphically isometrically into $C^k\times M'$ with respect to
the metric at time $t>0$, where
$M'$ is a complete simply connected K\"ahler manifold of nonnegative
bisectional curvature and positive Ricci curvature. For the compact cases,
there are results of this type in \cite{H-S-W},
see also \cite{B} and \cite{M1, p.64}.

2. For the ALE K\"ahler manifolds, a Liouville theorem on  plurisubharmonic
functions was proved earlier in \cite{N-S-T2}.
\endproclaim

\input amstex
\documentstyle{amsppt}
\magnification=1200
\hsize=13.8cm
\catcode`\@=11
\def\NoLogo{\let\logo@\empty}
\catcode`\@=\active
\NoLogo

\def\heat{\lf(\frac{\p}{\p t}-\Delta\ri)}
\def \b {\beta}
\def \d{\delta}
\def\i{\sqrt{-1}}

\def\lf{\left}
\def\ri{\right}
\def\bbar{\bar \beta}
\def\a{\alpha}
\def\ol{\overline}
\def\g{\gamma}
\def\e{\epsilon}
\def\p{\partial}
\def\delbar{{\bar\delta}}

\def\dbar{\bar\partial}

\def\C{\Bbb C}
\def\R{\Bbb R}

\def\vp{\varphi}

\def\dbar{\bar\partial}
\def\ba{{\bar\alpha}}
\def\bb{{\bar\beta}}
\def\bg{{\bar\gamma}}
\def\abb{{\alpha\bar\beta}}

\def\i{\sqrt {-1}}

\def \D {\Delta}
\def\aint{\frac{\ \ }{\ \ }{\hskip -0.4cm}\int}

\subheading{ \S3 Proof of  Theorem 2 and 3}

In order to prove Theorem 2
 we need a result in [L1, Corollary 5]. For the sake of completeness,
we will sketch the proof. It seems that in the proof of this result,
we need to assume that the holomorphic bisectional curvature is nonnegative.

\proclaim{Lemma 3.1} Let $M$ be a complete noncompact K\"ahler manifold with nonnegative holomorphic bisectional curvature. If $f$ is a harmonic function with sub-quadratic growth defined on $M$, then $f$ is pluri-harmonic.
\endproclaim
\demo{Proof}
Let $h=||f_{\abb}||^2=g^{\a\bar\delta}g^{\gamma\bar\beta}f_{\abb}f_{\gamma\bar\delta}$,  where $g_{\abb}$ is the metric of $M$ and $\g^{\abb}$ is its inverse.
Since $f$ is harmonic, by (1.8)
$$
\Delta f_{\gamma\bar\delta}=-R_{\beta\bar\a\gamma\bar\delta}f_{\abb}+\frac12\lf(R_{\gamma \bar p}f_{p\bar\delta}+R_{p\bar\delta}f_{\gamma\bar p}\ri).
$$
Hence in normal coordinates so that at a point $x$, $f_{\abb}=\lambda_{\a}\delta_{\a\b}$, we have
$$
\split
\Delta h & = 2f_{\gamma\bar\delta s\bar s}f_{\delta\bar \gamma}+||f_{\abb\gamma}||^2+||f_{\a\bb\bg}||^2\\
&=-2R_{\b\a\gamma\bar\delta}f_{\abb}f_{\delta\bar\gamma}+\lf(R_{\gamma \bar p}f_{p\bar\delta}+R_{p\bar\delta}f_{\gamma\bar p}\ri)f_{\delta\bar\gamma}+||f_{\abb\gamma}||^2+||f_{\a\bb\bg}||^2\\
&= -2R_{\a\ba\gamma\bar\gamma}\lambda_{\a}\lambda_{\gamma}+2R_{\gamma\bg}\lambda_\gamma^2+||f_{\abb\gamma}||^2+||f_{\a\bb\bg}||^2\\
&=\sum_{\a,\b}R_{\a\ba\b\bb}\lf(\lambda_{\a}-\lambda_{\b}\ri)^2+||f_{\abb\gamma}||^2+||f_{\a\bb\bg}||^2
\\
&\ge 0,
\endsplit
$$
where we have used the fact that $M$ has nonnegative holomorphic bisectional curvature. Since $|f(x)|=o\lf(r^2(x)\ri)$ where $r(x)$ is the distance from a fixed point $o\in M$, as in [L1, p.90-91], we have
$$
\frac{1}{V_o(R)}\int_{B_o(R)}h\le  \frac{C}{R^{-2}V_o(R)}\int_{B_o(R)}|\nabla f|^2=o(1),
$$
as $R\to\infty$. Here   $C$ is a constant independent of $R$ and we has used the gradient estimate in [C-Y]. Since $h$ is subharmonic, $h\equiv0$ by the mean value inequality in [L-S]. Hence $f$ is pluri-harmonic.
\enddemo

\proclaim{Lemma 3.2}
Let $M$ be a complete noncompact K\"ahler manifold with nonnegative holomorphic bisectional curvature.
Let $f$ be a pluri-harmonic function.
Then $\log (1+|\nabla f|^2)$ is pluri-subharmonic.
 \endproclaim
\demo{Proof} We adapt the complex notation. Let $h=|\nabla f|^2=\sum_{\a,\b}g^{\abb}f_{\a}f_{\bar{\b}}$. Here $g_{\abb}$ is the K\"ahler metric and $(g^{\abb})$ is the inverse of $(g_{\abb})$.
To prove that $\log (1+h)$ is pluri-subharmonic,
it is sufficient to  show that $\lf[\log (1+h)\ri]_{\g\bg}\ge 0$ in
 normal coordinates. Direct calculation shows that:
$$
\split
h_{\g\bar{\g}} & =  \left(\sum_{\a\b}g^{\abb}f_{\a}f_{\bar{\b}}\right)_{\g\bar{\g}}\\
& =  \sum_{\a,\b}g^{\abb}\lf[f_{\a\g}f_{\bar{\b}\bar{\g}}+f_{\a\bar\g}f_{\bar\b\g}+f_{\a\g\bar{\g}}f_{\bar{\b}}
+f_{\a}f_{\bar{\b}\g \bar{\g}}\ri]\\
& = \sum_{\a}  f_{\a\g}f_{\bar{\a}\bar{\g}}+ \sum_{\a,s}  R_{\g\bar{\g}\a \bar{s}}f_{s}f_{\bar{\a}}
\endsplit \tag 3.1
$$
where we have used the fact that $f$ is pluri-harmonic. Hence
$$
\split
\lf[\log (1+h)\ri]_{\g\bg}&=\frac{1}{(1+h)^2}
\lf[(1+h)h_{\g\bg}-h_\g h_\bg\ri]\\
&=\frac{1}{(1+h)^2}\bigg[(1+h)\lf(\sum_{\a}  f_{\a\g}f_{\bar{\a}\bar{\g}}+ \sum_{\a,s}  R_{\g\bar{\g}\a \bar{s}}f_{s}f_{\bar{\a}}\ri)\\
& \quad-\sum_\a  f_{\a\g}f_{\ba}\sum_\a f_\a f_{\ba\bg} \bigg]\\
&\ge\frac{1}{(1+h)^2} \lf(\sum_{\a}  f_{\a\g}f_{\bar{\a}\bar{\g}}+ \sum_{\a,s}  R_{\g\bar{\g}\a \bar{s}}f_{s}f_{\bar{\a}}\ri)
\endsplit
\tag 3.2
$$
where we have used the fact that $f$ is pluri-harmonic.
From (3.2), the fact that  $M$ has nonnegative
 holomorphic bisectional curvature, it is easy to see that $\log (1+h)$ is
pluri-subharmonic.
\enddemo

\proclaim{Lemma 3.3} Let $M$ be a complete noncompact K\"ahler manifold with nonnegative holomorphic bisectional curvature such that $||Rm||$ and $||\nabla Rm||$ are bounded.
Let $f$ be a harmonic function on $M$ satisfying (0.1).  Then
$$
||f_{\a\b}||(x)\le C\lf(1+r(x)\ri)^{3/2}\tag 3.3
$$
for some constant $C$, where $r(x)$ is the distance from $x$ to a fixed point $o\in M$.
\endproclaim

\demo{Proof} By Lemma 3.1, $f$ is pluri-harmonic. Let $\Psi =||f_{\a\b}||^2$. Then in normal coordinates:
$$
\split
\Delta \Psi & =  \sum_{s}\lf(\sum_{\a,\b,\g,\d } g^{\a \bar \g}g^{\b\bar\d}f_{\a\b}f_{\bar\g\bar\d}\ri)_{s\bar s} \\
& =  \sum_{s}\sum_{\a,\b,\g,\d } g^{\a \bar \g}g^{\b\bar\d}\lf(f_{\a\b s\bar s}f_{\bar\g\bar\d}+f_{\bar\g\bar\d s\bar s}f_{\a\b}+f_{\a\b s}f_{\bar\g\bar\d \bar s}+f_{\bar \g\bar \d s}f_{\a\b \bar s}\ri)\\
&\ge  \sum_{\a,\b, s}f_{\a\b s\bar s}f_{\bar \a\bar \b}+f_{\bar \a\bar\b s\bar s}f_{\a\b}\\
&\ge \sum_{\a,\b, s}\lf(R_{\a\bar{s}, \beta}f_{s}f_{\bar{\a}\bar{\b}}
+R_{s\bar{\a},\bar{\b}}f_{\bar{s}}f_{\a\b}\ri)+2\sum_{\a,\b, s,t} R_{\a\bar t\b\bar s}f_{ts}f_{\bar \a\bar\b}
\endsplit
$$
where $R_{\a\bar \b}$ is the Ricci tensor of $M$. Hence
$$
\Delta \Psi\ge -C_1(\Psi+h)
$$
for some constant $C_1>0$ depending only on $m$ and the bound of $||Rm||+||\nabla Rm||$. By ( 3.1), we also have
$$
\Delta h\ge \Psi.
$$
Let $S_R=\sup_{B_o(R)}h$, then

$$
\Delta(\Psi+C_1h)\ge -C_1S_{2R}
$$
on $B_o(2R)$. Hence for any $T>0$, we have
$$
\lf(\D-\frac{\p}{\p t}\ri)\lf(\Psi+C_1h+C_1S_{2R}(T-t)\ri)\ge0.
$$
Since $\Psi+C_1h+C_1S_{2R}(T-t)\ge0$ for $0\le t\le T$, by [L-T, Theorem
1.1], for any $R>0$, if we let   $T=\frac14 R^2$, we have
$$
\split
\sup_{B_o(\frac 12 R)\times[\frac 18 R^2,\frac 14 R^2]}& \lf(\Psi+C_1h+C_1S_{2R}(T-t)\ri)\\
 &\le \frac{C_2}{R^2V_o(R)}\int_{\frac1{16} R^2}^{\frac14 R^2}\int_{B_o(R)}\lf(\Psi+C_1h+C_1S_{2R}(T-t)\ri) dx dt\\
&\le C_3\lf(R^{-2}+S_{2R}R^2\ri)
\endsplit
$$
for some constants $C_2$, $C_3>0$ independent of $R$. Here we have used the fact that
$$
\frac{1}{V_o(R)}\int_{B_o(R)}\Psi\le CR^{-2}
$$
for some constant $C$ independent of $R$, see [L1, p.90-91].
Since $S_{2R}=o\lf( R^{1/2}\ri)$ by the gradient estimate in [C-Y],  we have
$$
\Psi(x)\le C_4\lf(1+r(x)\ri)^{3}.
$$
\enddemo

Now we are ready to prove Theorem 2. \vskip .2cm

\demo{Proof of Theorem 2}
Let $f$ be a nonconstant   harmonic function on $M$ satisfying (0.1).
Then $f$ is pluri-harmonic by Lemma 3.1.
Since $M$ has bounded curvature, we can solve the K\"ahler-Ricci flow
$$
\frac{\p \tilde g_{\abb}}{\p t}=-\tilde R_{\abb}
$$
where initial data $\tilde g_{\abb}(\cdot,0)=g_{\abb}$.
The equation has a short time solution so that for any fixed $t>0$,
$(M,\tilde g_{\abb})$ still has nonnegative holomorphic bisectional
curvature  so that the curvature tensor and the covariant derivatives of the curvature tensor are bounded. Moreover, $\tilde g_{\abb}$ is uniformly equivalent to $g_{\abb}$. All these results are in [Sh2].
Fix $t>0$, then $f$ is still  a pluri-harmonic function on $(M,\tilde g_{\abb})$ satisfying (0.1). By Lemma 3.2, the function
$u=\log(1+|\nabla^{(t)}f|^2)$ is pluri-subharmonic, where $\nabla^{(t)}$ is the gradient with respect to the metric $\tilde g_{\abb}(\cdot,t)$.
By the gradient estimates in [C-Y], $|u|(x)=o(\log r(x))$.
 By Lemma 3.3  and (3.2), we have
$$
\Delta^{(t)} u\le C\lf(1+r(x)\ri)^3
$$
if $r(x)>1$, where $\Delta^{(t)}$ is the Laplacian for the metric $\tilde g_{\abb}(\cdot,t)$. By Theorem 1,
we conclude that $|\nabla^{(t)} f|$ depends only on $t$, where $\nabla^{(t)}$ is the gradient with respect to $\tilde g_{\abb}(\cdot,t)$. Let $t\to0$, we conclude that $|\nabla f|$ is constant. Hence $f$ must be of linear growth. Moreover,  by the Bochner formula, we conclude that $\nabla f$ must be parallel.
Hence $J(\nabla f)$ is also parallel, where $J$ is the complex structure of $M$.  From this it is easy to see that the universal cover of $M$ splits as
$\widetilde M'\times\C$ isometrically and holomorphically. At the same time by integrating along $\nabla f$,
$M$ splits as $\widetilde M \times \R$ isometrically, where $\widetilde M$ can be taken  to be the  component of $f^{-1}(0)$. In this case that $M$ supports
a nonconstant holomorphic function of growth rate (0.1), both the real and
imaginary part will split a factor of $\R$ and clearly that they consist
a complex plane $\C$.

\enddemo

\demo{Proof of Theorem 3} By the assumptions that
$$
\int_0^r s\lf(\aint_{B_o(s)} \Cal R(y)dy\ri)ds=o\lf(\log r\ri), \tag 3.4
$$
and that $\Cal R(x)\le C\lf(r^2(x)+1\ri)$, it is easy to see that the conditions in   Theorem 5.1 of \cite{N-S-T1} are satisfied and so there exists a solution $u(x)$ to the Poincar\'e-Lelong  equation $\i\p\dbar u =Ric_M$ such that
$$
\limsup_{x\to\infty}\frac{u(x)}{ \log  r(x)}=0.
$$
 Obviously, $u$ is plurisubharmonic. By Theorem 1, $u$ must be constant and so $M$ must be flat.
\enddemo

Finally, we should point out that in order to solve the Poincar\'e-Lelong equation
we only need (3.4) together with
$\lim \inf_{r\to \infty}\aint_{B_o(r)}{\Cal R}^2(y)\, dy =0$ which is
slightly more general than the assumptions on ${\Cal R}(x)$ in Theorem 3. Therefore,
Theorem 3 holds under these more general assumptions.

\enddocument